\documentclass[11pt,a4paper]{article}    
\usepackage{a4,amssymb,latexsym,amsmath,color}
\begin{document}

% =======================================================================
% Macro
% =======================================================================

\newcounter{z}
\newenvironment{enum}{\setcounter{z}{0}
\begin{list}{\rm\alph{z})}{\usecounter{z}
\setlength{\topsep}{1ex}\setlength{\labelwidth}{0.5cm}
\setlength{\leftmargin}{1cm}\setlength{\labelsep}{0.25cm}
\setlength{\parsep}{-3pt}}}{\end{list}~\\[-6ex]}

\newcounter{z1}
\newenvironment{enum1}{\setcounter{z1}{0}
\begin{list}{\arabic{z1})}{\usecounter{z1}
\setlength{\topsep}{1ex}\setlength{\labelwidth}{0.6cm}
\setlength{\leftmargin}{1cm}\setlength{\labelsep}{0.25cm}
\setlength{\parsep}{-3pt}}}{\end{list}~\\[-6ex]}

\newcounter{z2}
\newenvironment{enum2}{\setcounter{z2}{0}
\begin{list}{$\bullet$}{\usecounter{z2}
\setlength{\topsep}{1ex}\setlength{\labelwidth}{0.6cm}
\setlength{\leftmargin}{0.75cm}\setlength{\labelsep}{0.25cm}
\setlength{\parsep}{-3pt}}}{\end{list}~\\[-6ex]}

\newcounter{zb}
\newenvironment{enumbib}{\setcounter{zb}{0}
\begin{small}
\begin{list}{\arabic{zb}.}{\usecounter{zb}
\setlength{\topsep}{1ex}\setlength{\labelwidth}{0.6cm}
\setlength{\leftmargin}{0.8cm}\setlength{\labelsep}{0.25cm}
\setlength{\parsep}{1pt}}}{\end{list}\end{small}}

\newcounter{pic}\setcounter{pic}{0}

\newcommand{\bdpm}{\begin{displaymath}}
\newcommand{\edpm}{\end{displaymath}}

\newcommand{\beas}{\begin{eqnarray*}}
\newcommand{\eeas}{\end{eqnarray*}}

\newcommand{\ba}{\begin{array}}
\newcommand{\ea}{\end{array}}

\newenvironment{proof}{\vspace*{-0.25cm}\begin{sloppypar}\noindent{\it 
Proof.}}{\hfill\qed\end{sloppypar}\vspace*{0.15cm}}

\newtheorem{theo}{Theorem}[section]      
\newtheorem{prop}[theo]{Proposition}
\newtheorem{lem}[theo]{Lemma}
\newtheorem{cor}[theo]{Corollary}
\newtheorem{exa}[theo]{Example}
\newtheorem{rem}[theo]{Remark}

\newcommand{\brm}{\begin{rm}}
\newcommand{\erm}{\end{rm}}

\newcommand{\qed}{\hfill $\Box$}
\newcommand{\sg}{{\sf sign}}
\newcommand{\inv}{{\sf inv}}
\newcommand{\lis}{{\sf lis}}
\newcommand{\lf}{\lfloor}
\newcommand{\rf}{\rfloor}
\newcommand{\nat}{{\Bbb N}}

% =======================================================================
% Title
% =======================================================================

\begin{center}
{\large\bf PERMUTATION SIGN UNDER THE}\\ 
{\large\bf ROBINSON-SCHENSTED-KNUTH CORRESPONDENCE}\\[1cm]
Astrid Reifegerste\\
Institut f\"ur Mathematik, Universit\"at Hannover\\
Welfengarten 1, 30167 Hannover, Germany\\
{\it reifegerste@math.uni-hannover.de}\\[0.5cm]
version of June 30, 2003 
\end{center}
\vspace*{0.5cm}

\begin{footnotesize}
{\sc Abstract.} We show how the sign of a permutation can be deduced from the 
tableaux induced by the permutation under the Robinson-Schensted-Knuth 
correspondence. The result yields a simple proof of a conjecture on the squares  
of imbalances raised by Stanley.
\end{footnotesize}
\vspace*{1cm}

% =======================================================================
% Section 1: Introduction
% =======================================================================

\setcounter{section}{1}\setcounter{theo}{0}

\centerline{\large{\bf 1}\hspace*{0.25cm}
{\sc Introduction}}
\vspace*{0.5cm}

The Robinson-Schensted-Knuth correspondence is a classical bijection between 
permutations and pairs of standard Young tableaux of the same shape. It was 
given in \cite{schensted}, and followed by numerous papers studying its 
combinatorial properties.\\
There are some well-known relations between certain permutation statistics on 
the one hand and tableaux statistics on the other hand. Schensted's classical 
Theorem \cite{schensted} states that the length of the longest increasing 
and decreasing subsequence in a permutation equals the length of the first row 
and column of the associated tableaux, respectively. Sch\"utzenberger 
\cite{schuetzenberger} showed that the descent set of a permutation and  
the descent set of its recording tableau coincide.\\[2ex] 
In \cite{reifegerste}, the question of how to obtain the sign from the 
associated pair of tableaux was answered for $321$-avoiding permutations, that 
is, for permutations having no decreasing subsequence of length three. The 
problem was motivated by the wish for refining the sign-balance property of $321$-avoiding 
permutations respecting the length of the longest increasing subsequence.\\ 
In this paper, we will give a general answer.\\[2ex]
It is known for a long time how the sign of an involution can be deduced from 
its corresponding tableau. A theorem of Sch\"utzenberger 
\cite{schuetzenberger1} says that the number of fixed points of an involution 
is equal to the number of columns of odd length in the associated tableau. 
Consequently, the sign is determined by the total length of all even-indexed tableau 
rows. In \cite{beissinger}, Beissinger described an elegant algorithm for 
constructing the tableau in bijection with an involution. Her bijective proof of Sch\"utzenberger's result makes use of the transparency 
of her algorithm with regard to the cycle structure.\\ 
Due to Knuth's equivalence relation, we can extend the case of involutions, treated in Section 3, to arbitrary permutations. To this end, 
we only have to consider the effect having an elementary Knuth transformation on the recording 
tableau. In Section 4, the resulting observation will be described, and the main result 
will be proved. Therefore, the sign of a permutation is the product of the signs of 
the associated tableaux and the parity of the total length of the even-indexed 
rows.\\[2ex]
In Section 5, we apply this result to prove Stanley's conjecture 
\cite{stanley1} on the squares of imbalances. Furthermore we use the explicit 
knowledge of the imbalance of hooks for determining the joint distribution of the 
sign and the length of the longest increasing subsequence for a special kind of 
pattern-avoiding permutations.  
\vspace*{0.75cm}

% =======================================================================
% Section 2: Definitions and notations
% =======================================================================

\setcounter{section}{2}\setcounter{theo}{0}

\centerline{\large{\bf 2}\hspace*{0.25cm}
{\sc Definitions and notations}}
\vspace*{0.5cm}

Let $\lambda=(\lambda_1,\lambda_2,\ldots,\lambda_l)$ with 
$\lambda_1\ge\lambda_2\ge\ldots\ge\lambda_l>0$ be a partition of $n$, denoted 
by $\lambda\vdash n$. We use this notation for both the partition and its 
corresponding Young diagram. A standard Young tableau (or just tableau in the 
following) of shape $\lambda$ is a labeling of the squares of $\lambda$ with 
the numbers $1,2,\ldots,n$ so that the rows and columns are increasing.\\[2ex]
We assume that the reader is familar with the combinatorics of Young tableaux, 
and, especially, with the Robinson-Schensted-Knuth correspondence (for details see, e.g., \cite{fulton} or 
\cite{stanley}).\\
This correspondence gives a bijection between permutations of the symmetric group ${\cal S}_n$ 
and pairs $(P,Q)$ of tableaux of the same shape $\lambda\vdash n$. Here $P$ is called the {\it insertion tableau} 
and $Q$ the {\it recording tableau} of the permutation. If necessary, we will denote the tableaux 
that are obtained at the $i$th step of the algorithm by $P_i$ and $Q_i$. (By 
definition, we have $P=P_n$ and $Q=Q_n$.)\\[2ex] 
For a tableau $T$, we call a pair $(i,j)$ of entries an {\it inversion of $T$} if $j<i$ 
and $j$ is contained in a row below the row of $i$. We write $\inv(T)$ to 
denote the number of inversions of $T$. Furthermore, we define the {\it sign of a 
tableau $T$} as 
\bdpm
\sg(T)=(-1)^{\inv(T)}.
\edpm
Any tableau $T$ can be identified with its {\it row word}, denoted by $w(T)$, which is obtained 
by reading the entries row by row from left to right, starting from the top. 
(Note that the inversions of a tableau are just the inversions of this word.) 
\newpage

% =======================================================================
% Section 3: Involutions
% =======================================================================

\setcounter{section}{3}\setcounter{theo}{0}

\centerline{\large{\bf 3}\hspace*{0.25cm}
{\sc The case of involutions -- Beissinger's algorithm}}
\vspace*{0.5cm}

An important property of the Robinson-Schensted-Knuth correspondence is its symmetry: the recording 
tableau for a permutation $\pi$ is just the insertion tableau for the inverse 
permutation $\pi^{-1}$. This statement has already appeared in the work of 
Robinson \cite{robinson}; the first proof was given by Sch\"utzenberger \cite{schuetzenberger}. Later, the result was demonstrated again by Viennot \cite{viennot} using the symmetry of his geometric construction. 
Thus each involution $\pi\in{\cal S}_n$ is associated with a pair $(P,P)$ where $P$ is a tableau of shape $\lambda\vdash n$. 
Consequently, the correspondence gives a bijection between involutions and tableaux.\\[2ex] 
The question of how to see the sign of an involution by looking at its 
corresponding tableau can be easily answered using a simple algorithm given 
by Beissinger \cite{beissinger}. Her construction yields the tableau in 
bijection with an involution, and works as follows.\\  
Given an involution $\pi\in{\cal S}_n$, write $\pi$ as a product of disjoint cycles $c_1,c_2,\ldots,c_k$ in increasing order of 
their greatest element and with $a<b$ for each cycle $(a,b)$. Starting with the 
empty tableau, $P$ is obtained recursively by applying the 
following procedure, for $i=1,\ldots,k$: 
\begin{enum2}
\item If $c_i=(a)$, then place $a$ at the end of the first row of $P$.
\item If $c_i=(a,b)$, then insert $a$ into $P$. Let $l$ be the index of the row in which the insertion 
process stops. Place $b$ at the end of the row $l+1$.
\end{enum2}
 
One consequence of the algorithm is a direct proof of a result of 
Sch\"utzenberger \cite[p. 93]{schuetzenberger1} which we give here in the following equivalent formulation.

\begin{theo} \label{invo}
Let $\pi\in{\cal S}_n$ be an involution. Then $\sg(\pi)=(-1)^e$ where $e$ 
denotes the total length of all even-indexed rows of the corresponding tableau.   
\end{theo}

\begin{proof}
By description of the algorithm, inserting the elements of a 2-cycle 
creates two squares in consecutive rows. Hence the number of 2-cycles of $\pi$ 
equals the number of squares in all even-indexed rows of $P$. This yields the 
assertion.
\end{proof}

\begin{exa}
\brm
Let $\pi=4\:8\:7\:1\:9\:6\:3\:2\:5=(1,4)(6)(3,7)(2,8)(5,9)\in{\cal S}_9$. 
Beissinger's algorithm generates the corresponding tableau in five steps:
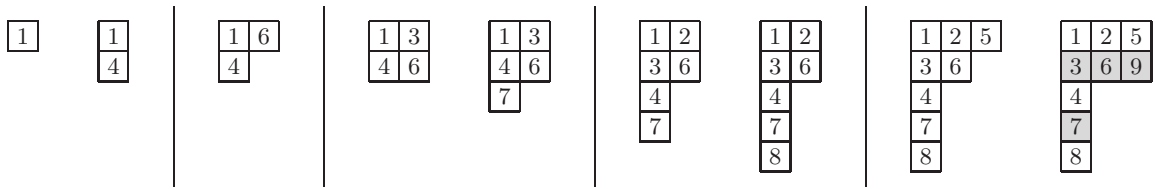
\begin{center}
% == Picture ==========================================
\unitlength=0.4cm
\definecolor{gray1}{gray}{0.85}
\fboxsep0cm
\fboxrule0cm
\begin{picture}(38,6)
\put(35,3.5){\fcolorbox{gray1}{gray1}{\makebox(3,1){}}}
\put(35,1.5){\fcolorbox{gray1}{gray1}{\makebox(1,1){}}}
\linethickness{0.45pt}
\multiput(0,4.5)(0,1){2}{\line(1,0){1}}
\multiput(0,4.5)(1,0){2}{\line(0,1){1}}
\multiput(3,3.5)(0,1){3}{\line(1,0){1}}
\multiput(3,3.5)(1,0){2}{\line(0,1){2}}
\multiput(7,4.5)(0,1){2}{\line(1,0){2}}
\multiput(7,3.5)(1,0){2}{\line(0,1){2}}
\put(7,3.5){\line(1,0){1}}
\put(9,4.5){\line(0,1){1}}
\multiput(12,3.5)(0,1){3}{\line(1,0){2}}
\multiput(12,3.5)(1,0){3}{\line(0,1){2}}
\multiput(16,3.5)(0,1){3}{\line(1,0){2}}
\multiput(16,2.5)(1,0){2}{\line(0,1){3}}
\put(16,2.5){\line(1,0){1}}
\put(18,3.5){\line(0,1){2}}
\multiput(21,3.5)(0,1){3}{\line(1,0){2}}
\multiput(21,1.5)(0,1){2}{\line(1,0){1}}
\multiput(21,1.5)(1,0){2}{\line(0,1){4}}
\put(23,3.5){\line(0,1){2}}
\multiput(25,3.5)(0,1){3}{\line(1,0){2}}
\multiput(25,0.5)(0,1){3}{\line(1,0){1}}
\multiput(25,0.5)(1,0){2}{\line(0,1){5}}
\put(27,3.5){\line(0,1){2}}
\multiput(30,4.5)(0,1){2}{\line(1,0){3}}
\put(30,3.5){\line(1,0){2}}
\multiput(30,0.5)(0,1){3}{\line(1,0){1}}
\multiput(30,0.5)(1,0){2}{\line(0,1){5}}
\put(32,3.5){\line(0,1){2}}
\put(33,4.5){\line(0,1){1}}
\multiput(35,3.5)(0,1){3}{\line(1,0){3}}
\multiput(35,0.5)(0,1){3}{\line(1,0){1}}
\multiput(35,0.5)(1,0){2}{\line(0,1){5}}
\multiput(37,3.5)(1,0){2}{\line(0,1){2}}
\put(0.5,5){\makebox(0,0)[cc]{\footnotesize 1}}
\put(3.5,5){\makebox(0,0)[cc]{\footnotesize 1}}
\put(3.5,4){\makebox(0,0)[cc]{\footnotesize 4}}
\put(7.5,5){\makebox(0,0)[cc]{\footnotesize 1}}
\put(8.5,5){\makebox(0,0)[cc]{\footnotesize 6}}
\put(7.5,4){\makebox(0,0)[cc]{\footnotesize 4}}
\put(12.5,5){\makebox(0,0)[cc]{\footnotesize 1}}
\put(13.5,5){\makebox(0,0)[cc]{\footnotesize 3}}
\put(12.5,4){\makebox(0,0)[cc]{\footnotesize 4}}
\put(13.5,4){\makebox(0,0)[cc]{\footnotesize 6}}
\put(16.5,5){\makebox(0,0)[cc]{\footnotesize 1}}
\put(17.5,5){\makebox(0,0)[cc]{\footnotesize 3}}
\put(16.5,4){\makebox(0,0)[cc]{\footnotesize 4}}
\put(17.5,4){\makebox(0,0)[cc]{\footnotesize 6}}
\put(16.5,3){\makebox(0,0)[cc]{\footnotesize 7}}
\put(21.5,5){\makebox(0,0)[cc]{\footnotesize 1}}
\put(22.5,5){\makebox(0,0)[cc]{\footnotesize 2}}
\put(21.5,4){\makebox(0,0)[cc]{\footnotesize 3}}
\put(22.5,4){\makebox(0,0)[cc]{\footnotesize 6}}
\put(21.5,3){\makebox(0,0)[cc]{\footnotesize 4}}
\put(21.5,2){\makebox(0,0)[cc]{\footnotesize 7}}
\put(25.5,5){\makebox(0,0)[cc]{\footnotesize 1}}
\put(26.5,5){\makebox(0,0)[cc]{\footnotesize 2}}
\put(25.5,4){\makebox(0,0)[cc]{\footnotesize 3}}
\put(26.5,4){\makebox(0,0)[cc]{\footnotesize 6}}
\put(25.5,3){\makebox(0,0)[cc]{\footnotesize 4}}
\put(25.5,2){\makebox(0,0)[cc]{\footnotesize 7}}
\put(25.5,1){\makebox(0,0)[cc]{\footnotesize 8}}
\put(30.5,5){\makebox(0,0)[cc]{\footnotesize 1}}
\put(31.5,5){\makebox(0,0)[cc]{\footnotesize 2}}
\put(32.5,5){\makebox(0,0)[cc]{\footnotesize 5}}
\put(30.5,4){\makebox(0,0)[cc]{\footnotesize 3}}
\put(31.5,4){\makebox(0,0)[cc]{\footnotesize 6}}
\put(30.5,3){\makebox(0,0)[cc]{\footnotesize 4}}
\put(30.5,2){\makebox(0,0)[cc]{\footnotesize 7}}
\put(30.5,1){\makebox(0,0)[cc]{\footnotesize 8}}
\put(35.5,5){\makebox(0,0)[cc]{\footnotesize 1}}
\put(36.5,5){\makebox(0,0)[cc]{\footnotesize 2}}
\put(37.5,5){\makebox(0,0)[cc]{\footnotesize 5}}
\put(35.5,4){\makebox(0,0)[cc]{\footnotesize 3}}
\put(36.5,4){\makebox(0,0)[cc]{\footnotesize 6}}
\put(37.5,4){\makebox(0,0)[cc]{\footnotesize 9}}
\put(35.5,3){\makebox(0,0)[cc]{\footnotesize 4}}
\put(35.5,2){\makebox(0,0)[cc]{\footnotesize 7}}
\put(35.5,1){\makebox(0,0)[cc]{\footnotesize 8}}
\linethickness{0.3pt}
\put(5.5,0){\line(0,1){6}}
\put(10.5,0){\line(0,1){6}}
\put(19.5,0){\line(0,1){6}}
\put(28.5,0){\line(0,1){6}}
\end{picture}
\end{center}
Its even-indexed rows consist of 4 squares in all. Thus $\pi$ is even.
\erm
\end{exa}

\begin{rem}
\brm
Another proof of Theorem \ref{invo} is given in \cite[Exercise 7.28a]{stanley}. 
It is based on growth diagrams of permutations.
\erm
\end{rem}
\vspace*{0.5cm}

% =======================================================================
% Section 4: General case
% =======================================================================

\setcounter{section}{4}\setcounter{theo}{0}

\centerline{\large{\bf 4}\hspace*{0.25cm}
{\sc The general case -- Knuth equivalence}}
\vspace*{0.5cm}

Now we extend the result given for involutions in the previous section to 
arbitrary permutations. To this end, another classical property of the 
Robinson-Schensted-Knuth correspondence will be used. In \cite{knuth}, Knuth gave
a characterization of permutations having the same insertion tableau.\\[2ex]   
An {\it elementary Knuth transformation} of a permutation is the application of 
one of the transformations 
\bdpm
acb\mapsto cab,\quad cab\mapsto acb,\quad bac\mapsto bca,\quad bca\mapsto 
bac\qquad\mbox{where }a<b<c
\edpm
to three consecutive letters of the permutation. Two permutations are called {\it 
Knuth-equivalent} if they can be changed into each other by a sequence of 
elementary Knuth transformations.\\[2ex]
The result of Knuth \cite[Theorem 6]{knuth} is the following: permutations are Knuth-equivalent if and only if their insertion tableaux 
coincide. Because of the symmetry of the Robinson-Schensted-Knuth correspondence, this also implies a 
characterization of permutations having the same recording tableau.\\[2ex] 
The proof of the generalization of Theorem \ref{invo} bases on the following 
simple observation.

\begin{lem} \label{knuth}
Any elementary Knuth transformation reverses the sign of the recording tableau.
\end{lem}

\begin{proof}
Let $\pi\in{\cal S}_n$ be a permutation and $(P,Q)$ the corresponding pair of 
tableaux. We show that the recording tableau $Q'$ of a permutation 
$\sigma$ which is obtained from $\pi$ by an elementary Knuth transformation arises from 
$Q$ by exchanging two consecutive numbers.\\ 
By symmetry, $Q$ and $Q'$ are the insertion tableaux of $p=\pi^{-1}$ and $s=\sigma^{-1}$, 
respectively. First, let 
\bdpm
s=p_1\cdots p_{a-1}\:k+2\:p_{a+1}\cdots p_{b-1}\:k\:p_{b+1}\cdots p_{c-1}\:k+1\:
p_{c+1}\cdots p_n=(k+1,k+2)\cdot p
\edpm
where $a<b<c$. Then $\sigma$ has arose from $\pi$ by the transformation 
$bac\mapsto bca$. Clearly, we have $Q_i=Q'_i$ for $i=1,\ldots,a-1$. Of course, the 
smallest integer in the 
first row of $Q_{a-1}$ that is larger than $k+1$ is just as the smallest integer that is larger than $k+2$. Thus the tableaux $Q_a=Q_{a-1}\leftarrow k+1$ and 
$Q'_a=Q_{a-1}\leftarrow k+2$ are of the same shape. 
(In the following, we use the notation $T\leftarrow i$ as well for the tableau resulting from the row insertion of the 
integer $i$ into the tableau $T$ as for the insertion process itself.)
More exactly, $Q_a$ and $Q'_a$ only differ in the 
entry of a single square (in the first row) which is occupied by $k+1$ in $Q_a$ and by 
$k+2$ in $Q'_a$. Since $p_i\notin\{k+1,k+2\}$ for $i=a+1,\ldots,c-1$, the insertion 
$Q_{i-1}\leftarrow p_i$ takes the same path as $Q'_{i-1}\leftarrow p_i$. 
(If $k+1$ is the smallest integer in $Q_{i-1}$ which is larger than 
an element $p_i$, then $k+2$ is it as well in $Q'_{i-1}$.) In this way, we obtain the tableaux
\begin{center}
% == Picture ==========================================
\unitlength=0.5cm
\definecolor{gray1}{gray}{0.85}
\fboxsep0cm
\fboxrule0cm
\begin{picture}(21,4)
\linethickness{0.5pt}
\put(0.25,2){\makebox(0,0)[lc]{$Q_{c-1}\;=$}}
\put(4,4){\line(1,0){5}}\put(9,3){\line(0,1){1}}
\put(7,3){\line(1,0){2}}\put(7,1){\line(0,1){2}}
\put(6,1){\line(1,0){1}}\put(6,0){\line(0,1){1}}
\put(4,0){\line(1,0){2}}\put(4,0){\line(0,1){4}}
\put(5,1.5){\fcolorbox{gray1}{gray1}{\makebox(0.5,0.5){}}}
\put(12.25,2){\makebox(0,0)[lc]{$Q'_{c-1}\;=$}}
\put(16,4){\line(1,0){5}}\put(21,3){\line(0,1){1}}
\put(19,3){\line(1,0){2}}\put(19,1){\line(0,1){2}}
\put(18,1){\line(1,0){1}}\put(18,0){\line(0,1){1}}
\put(16,0){\line(1,0){2}}\put(16,0){\line(0,1){4}}
\put(17,1.5){\fcolorbox{gray1}{gray1}{\makebox(0.5,0.5){}}}
\put(7,0){\makebox(0,0)[lc]{\footnotesize$k+1$}}
\put(19,0){\makebox(0,0)[lc]{\footnotesize$k+2$}}
\linethickness{0.3pt}
\put(17.25,1.75){\circle*{0.15}}
\put(5.25,1.75){\circle*{0.15}}
\put(5.25,1.75){\line(1,-1){1.5}}
\put(17.25,1.75){\line(1,-1){1.5}}
\end{picture}
\end{center}
which are identical, except from the entry of the gray square. Note that this 
square cannot belong to the first row since $k+1$ and $k+2$, respectively, have been bumped 
at the latest when inserting $p_b=k$. Furthermore, the entry $k$ is placed 
above $k+1$ and $k+2$, respectively. If $k$ is bumped into the row containing 
$k+1$ or $k+2$, then it for its part bumps these numbers. Therefore, the insertions 
$Q_{c-1}\leftarrow k+2$ and $Q'_{c-1}\leftarrow k+1$ run identically. 
Finally, inserting the remaining elements $p_{c+1},\ldots,p_n$ successively 
yields the tableaux $Q$ and $Q'$ which can be tranformed into each other by changing the 
elements $k+1$ and $k+2$: 
\begin{center}
% == Picture ==========================================
\unitlength=0.5cm
\definecolor{gray1}{gray}{0.85}
\fboxsep0cm
\fboxrule0cm
\begin{picture}(22,5)
\linethickness{0.5pt}
\put(0,2.5){\makebox(0,0)[lc]{$Q\;=$}}
\put(3,5){\line(1,0){5}}\put(8,3){\line(0,1){2}}
\put(6,3){\line(1,0){2}}\put(6,2){\line(0,1){1}}
\put(5,2){\line(1,0){1}}\put(5,1){\line(0,1){1}}
\put(4,1){\line(1,0){1}}\put(4,0){\line(0,1){1}}
\put(3,0){\line(1,0){1}}\put(3,0){\line(0,1){5}}
\put(4,1.5){\fcolorbox{gray1}{gray1}{\makebox(0.5,0.5){}}}
\put(6.5,4){\fcolorbox{gray1}{gray1}{\makebox(0.5,0.5){}}}
\put(13,2.5){\makebox(0,0)[lc]{$Q'\;=$}}
\put(16,5){\line(1,0){5}}\put(21,3){\line(0,1){2}}
\put(19,3){\line(1,0){2}}\put(19,2){\line(0,1){1}}
\put(18,2){\line(1,0){1}}\put(18,1){\line(0,1){1}}
\put(17,1){\line(1,0){1}}\put(17,0){\line(0,1){1}}
\put(16,0){\line(1,0){1}}\put(16,0){\line(0,1){5}}
\put(17,1.5){\fcolorbox{gray1}{gray1}{\makebox(0.5,0.5){}}}
\put(19.5,4){\fcolorbox{gray1}{gray1}{\makebox(0.5,0.5){}}}
\put(5.7,0){\makebox(0,0)[lc]{\footnotesize$k+1$}}
\put(8.4,2.25){\makebox(0,0)[lc]{\footnotesize$k+2$}}
\put(18.7,0){\makebox(0,0)[lc]{\footnotesize$k+2$}}
\put(21.4,2.25){\makebox(0,0)[lc]{\footnotesize$k+1$}}
\linethickness{0.3pt}
\put(4.25,1.75){\line(1,-1){1.25}}
\put(6.75,4.25){\line(1,-1){1.5}}
\put(17.25,1.75){\line(1,-1){1.25}}
\put(19.75,4.25){\line(1,-1){1.5}}
\put(4.25,1.75){\circle*{0.15}}
\put(6.75,4.25){\circle*{0.15}}
\put(17.25,1.75){\circle*{0.15}}
\put(19.75,4.25){\circle*{0.15}}
\end{picture}
\end{center}
Why? Since $k$ has been inserted before $k+2$ in $Q$ and 
before $k+1$ in $Q'$, it is not possible that $k+2$ and $k+1$ are bumped into 
the row which contains $k+1$ and $k+2$, respectively. When inserting 
$p_{c+1},\ldots,p_n$ successively into $Q_c$, the element $k+2$ can be bumped 
into the row containing $k$ at most. Then it occupies the square directly to 
the right of that of $k$. By construction, $k+2$ does not move while it is contained in 
the same row as $k$. Analogously, $k+1$ cannot pass by $k$ in $Q'$. Consequently, 
we have $\inv(Q)=\inv(Q')+1$.\\
Consider now the second kind of an elementary Knuth transformation: 
$acb\mapsto cab$. Let 
\bdpm
s=p_1\cdots p_{a-1}\:k+1\:p_{a+1}\cdots p_{b-1}\:k+2\:p_{b+1}\cdots p_{c-1}\:k\:
p_{c+1}\cdots p_n=(k,k+1)\cdot p
\edpm
for some $a<b<c$. By reasoning similarly as above, we can show that
\begin{center}
% == Picture ==========================================
\unitlength=0.5cm
\definecolor{gray1}{gray}{0.85}
\fboxsep0cm
\fboxrule0cm
\begin{picture}(21,4)
\linethickness{0.5pt}
\put(5.5,3.5){\fcolorbox{gray1}{gray1}{\makebox(0.5,0.5){}}}
\put(4.5,2.5){\fcolorbox{gray1}{gray1}{\makebox(0.5,0.5){}}}
\put(3.5,1){\fcolorbox{gray1}{gray1}{\makebox(0.5,0.5){}}}
\put(0,2){\makebox(0,0)[lc]{$Q_c\;=$}}
\put(3,4){\line(1,0){5}}\put(8,3){\line(0,1){1}}
\put(7,3){\line(1,0){1}}\put(7,2){\line(0,1){1}}
\put(5,2){\line(1,0){2}}\put(5,1){\line(0,1){1}}
\put(4,1){\line(1,0){1}}\put(4,0){\line(0,1){1}}
\put(3,0){\line(1,0){1}}\put(3,0){\line(0,1){4}}
\put(4.75,2.75){\circle*{0.15}}
\put(5.75,3.75){\circle*{0.15}}
\put(3.75,1.25){\circle*{0.15}}
\put(4.7,0.25){\makebox(0,0)[lc]{\footnotesize$k$}}
\put(5.5,1.25){\makebox(0,0)[lc]{\footnotesize$k+2$}}
\put(7.5,2.1){\makebox(0,0)[lc]{\footnotesize$k+1$}}
\put(4.75,2.75){\line(1,-1){1}}
\put(5.75,3.75){\line(1,-1){1.5}}
\put(3.75,1.25){\line(1,-1){0.75}}
\put(17.5,3.5){\fcolorbox{gray1}{gray1}{\makebox(0.5,0.5){}}}
\put(16.5,2.5){\fcolorbox{gray1}{gray1}{\makebox(0.5,0.5){}}}
\put(15.5,1){\fcolorbox{gray1}{gray1}{\makebox(0.5,0.5){}}}
\put(12,2){\makebox(0,0)[lc]{$Q'_c\;=$}}
\put(15,4){\line(1,0){5}}\put(20,3){\line(0,1){1}}
\put(19,3){\line(1,0){1}}\put(19,2){\line(0,1){1}}
\put(17,2){\line(1,0){2}}\put(17,1){\line(0,1){1}}
\put(16,1){\line(1,0){1}}\put(16,0){\line(0,1){1}}
\put(15,0){\line(1,0){1}}\put(15,0){\line(0,1){4}}
\put(16.75,2.75){\circle*{0.15}}
\put(17.75,3.75){\circle*{0.15}}
\put(15.75,1.25){\circle*{0.15}}
\put(16.7,0.25){\makebox(0,0)[lc]{\footnotesize$k+1$}}
\put(17.5,1.25){\makebox(0,0)[lc]{\footnotesize$k+2$}}
\put(19.5,2.1){\makebox(0,0)[lc]{\footnotesize$k$}}
\put(16.75,2.75){\line(1,-1){1}}
\put(17.75,3.75){\line(1,-1){1.5}}
\put(15.75,1.25){\line(1,-1){0.75}}
\end{picture}
\end{center}
All the numbers $i\not=k,k+1$ are placed at the same position in $Q'_c$ as they have in $Q_c$. Note that $k+2$ can be bumped into the row of 
$k$ in $Q_c$ and $k+1$ in $Q'_c$ at most, respectively; in this case, the 
numbers occupy adjacent squares. Let us consider what happens when 
we insert the remaining elements into $Q_c$ and $Q'_c$. If $k+2$ is above $k$ 
in $Q_{i-1}$ and $k+1$ in $Q'_{i-1}$, respectively, then $k+1$ bumps $k+2$ 
during the insertion $Q_{i-1}\leftarrow p_i$ ($i>c$) if and only if $k$ bumps 
$k+2$ during $Q'_{i-1}\leftarrow p_i$. Thus we have $\inv(Q)=\inv(Q')+1$ if the 
numbers $k$, $k+1$, and $k+2$ have the relative positions in $Q$ and $Q'$ as shown 
in the above figure. Assume now that $k+2$ has been bumped into the row of $k$ and $k+1$, 
respectively:
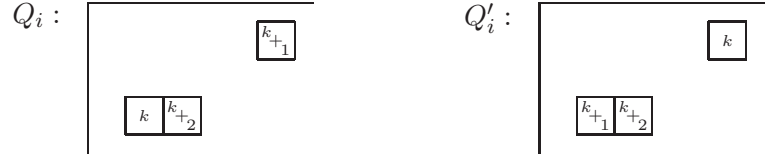
\begin{center}
% == Picture ==========================================
\unitlength=0.5cm
\begin{picture}(20,4)
\linethickness{0.5pt}
\put(0,4){\makebox(0,0)[lt]{$Q_i:$}}
\put(2,4){\line(1,0){6}}\put(2,0){\line(0,1){4}}
\multiput(3,0.5)(0,1){2}{\line(1,0){2}}
\multiput(3,0.5)(1,0){3}{\line(0,1){1}}
\multiput(6.5,2.5)(0,1){2}{\line(1,0){1}}
\multiput(6.5,2.5)(1,0){2}{\line(0,1){1}}
\put(3.5,1){\makebox(0,0)[cc]{\tiny$k$}}
\put(4.5,1){\makebox(0,0)[cc]{\tiny$+$}}
\put(4.25,1.25){\makebox(0,0)[cc]{\tiny$k$}}
\put(4.75,0.75){\makebox(0,0)[cc]{\tiny$2$}}
\put(7,3){\makebox(0,0)[cc]{\tiny$+$}}
\put(6.75,3.25){\makebox(0,0)[cc]{\tiny$k$}}
\put(7.25,2.75){\makebox(0,0)[cc]{\tiny$1$}}
\put(12,4){\makebox(0,0)[lt]{$Q'_i:$}}
\put(14,4){\line(1,0){6}}\put(14,0){\line(0,1){4}}
\multiput(15,0.5)(0,1){2}{\line(1,0){2}}
\multiput(15,0.5)(1,0){3}{\line(0,1){1}}
\multiput(18.5,2.5)(0,1){2}{\line(1,0){1}}
\multiput(18.5,2.5)(1,0){2}{\line(0,1){1}}
\put(15.25,1.25){\makebox(0,0)[cc]{\tiny$k$}}
\put(15.5,1){\makebox(0,0)[cc]{\tiny$+$}}
\put(15.75,0.75){\makebox(0,0)[cc]{\tiny$1$}}
\put(16.25,1.25){\makebox(0,0)[cc]{\tiny$k$}}
\put(16.5,1){\makebox(0,0)[cc]{\tiny$+$}}
\put(16.75,0.75){\makebox(0,0)[cc]{\tiny$2$}}
\put(19,3){\makebox(0,0)[cc]{\tiny$k$}}
\end{picture}
\end{center}
If $k+1$ and $k$ is bumped into this row by a sequence of insertions, 
respectively, then we obtain
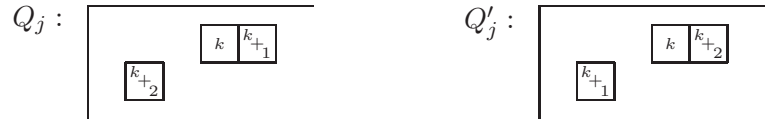
\begin{center}
% == Picture ==========================================
\unitlength=0.5cm
\begin{picture}(20,3)
\linethickness{0.5pt}
\put(0,3){\makebox(0,0)[lt]{$Q_j:$}}
\put(2,3){\line(1,0){6}}\put(2,0){\line(0,1){3}}
\multiput(5,1.5)(0,1){2}{\line(1,0){2}}
\multiput(5,1.5)(1,0){3}{\line(0,1){1}}
\multiput(3,0.5)(0,1){2}{\line(1,0){1}}
\multiput(3,0.5)(1,0){2}{\line(0,1){1}}
\put(5.5,2){\makebox(0,0)[cc]{\tiny$k$}}
\put(3.5,1){\makebox(0,0)[cc]{\tiny$+$}}
\put(3.25,1.25){\makebox(0,0)[cc]{\tiny$k$}}
\put(3.75,0.75){\makebox(0,0)[cc]{\tiny$2$}}
\put(6.5,2){\makebox(0,0)[cc]{\tiny$+$}}
\put(6.25,2.25){\makebox(0,0)[cc]{\tiny$k$}}
\put(6.75,1.75){\makebox(0,0)[cc]{\tiny$1$}}
\put(12,3){\makebox(0,0)[lt]{$Q'_j:$}}
\put(14,3){\line(1,0){6}}\put(14,0){\line(0,1){3}}
\multiput(17,1.5)(0,1){2}{\line(1,0){2}}
\multiput(17,1.5)(1,0){3}{\line(0,1){1}}
\multiput(15,0.5)(0,1){2}{\line(1,0){1}}
\multiput(15,0.5)(1,0){2}{\line(0,1){1}}
\put(15.25,1.25){\makebox(0,0)[cc]{\tiny$k$}}
\put(15.5,1){\makebox(0,0)[cc]{\tiny$+$}}
\put(15.75,0.75){\makebox(0,0)[cc]{\tiny$1$}}
\put(18.25,2.25){\makebox(0,0)[cc]{\tiny$k$}}
\put(18.5,2){\makebox(0,0)[cc]{\tiny$+$}}
\put(18.75,1.75){\makebox(0,0)[cc]{\tiny$2$}}
\put(17.5,2){\makebox(0,0)[cc]{\tiny$k$}}
\end{picture}
\end{center}
Clearly, the bumped numbers occur at the same position in the next row since 
they are consecutive. Each further insertion will take the same path through 
$Q_j$ as through $Q'_j$. Note again that $k+1$ and $k+2$, respectively, does 
not move while $k$ is situated to their left. Thus $k+1$ and $k+2$ are 
contained in different rows, and we have $\inv(Q)=\inv(Q')-1$.
\end{proof}

\begin{rem}
\brm
Obviously, the inversion numbers of permutations which are connected by an 
elementary Knuth relation differ by $1$. Since a transformation 
$bac(b)\mapsto bca(b)$ can cause as well the increase as the decrease (by $1$) 
of the inversion number of the recording tableau it seems that the problem of 
recovering the inversion number from the tableaux is unlike more difficult.  
\erm
\end{rem}

\begin{theo} \label{main}
Let $\pi\in{\cal S}_n$ be a permutation and $(P,Q)$ its associated pair of 
tableaux. Then
\bdpm
\sg(\pi)=\sg(P)\cdot\sg(Q)\cdot(-1)^e
\edpm 
where $e$ is the total length of all even-indexed rows of $P$.
\end{theo}

\begin{proof}
Consider two special elements of the Knuth class containing $\pi$: the 
involution $\sigma$, and the permutation $\tau$ whose recording tableau, 
denoted by $I$, has no inversions. (Note that the row word of $I$ is the 
identity in ${\cal S}_n$.)
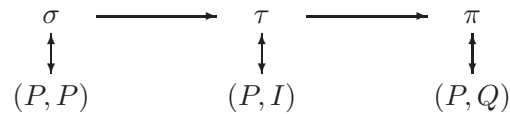
\begin{center}
% == Picture ==========================================
\unitlength=0.4cm
\begin{picture}(15,3)
\put(0.5,2.5){\makebox(0,0)[cb]{$\sigma$}}
\put(0.5,0){\makebox(0,0)[cc]{$(P,P)$}}
\put(7.5,2.5){\makebox(0,0)[cb]{$\tau$}}
\put(7.5,0){\makebox(0,0)[cc]{$(P,I)$}}
\put(14.5,2.5){\makebox(0,0)[cb]{$\pi$}}
\put(14.5,0){\makebox(0,0)[cc]{$(P,Q)$}}
\put(2,2.7){\vector(1,0){4}}
\put(9,2.7){\vector(1,0){4}}
\put(0.5,1.5){\vector(0,1){0.65}}
\put(0.5,1.5){\vector(0,-1){0.7}}
\put(7.5,1.5){\vector(0,1){0.65}}
\put(7.5,1.5){\vector(0,-1){0.7}}
\put(14.5,1.5){\vector(0,1){0.65}}
\put(14.5,1.5){\vector(0,-1){0.7}}
\end{picture}
\end{center}
Obviously, the tableau $P$ can be transformed into $I$ by $\inv(P)$ exchanges 
of consecutive elements $i$ and $i+1$ (necessarily, contained in different rows). In terms of row words, each of these 
transformations means the multiplication of $w(P)$ with an adjacent 
transposition $(i,i+1)$ in ${\cal S}_n$. Analogously, $Q$ arises from $I$ by 
$\inv(Q)$ exchanges of consecutive elements. Therefore, the assertion follows 
immediately from Theorem \ref{invo} and Lemma \ref{knuth}.
\end{proof}

\begin{exa}
\brm
Consider the permutation $\pi=2\:9\:1\:5\:6\:4\:8\:3\:7\in{\cal S}_9$. It 
corresponds to the pair 
\begin{center}
% == Picture ==========================================
\unitlength=0.4cm
\begin{picture}(16,4)
\linethickness{0.45pt}
\put(0,2){\makebox(0,0)[lc]{$P\;=$}}
\multiput(3,3)(0,1){2}{\line(1,0){4}}
\put(3,2){\line(1,0){3}}
\multiput(3,0)(1,0){2}{\line(0,1){4}}
\multiput(3,0)(0,1){2}{\line(1,0){1}}
\multiput(5,2)(1,0){2}{\line(0,1){2}}
\put(7,3){\line(0,1){1}}
\put(10,2){\makebox(0,0)[lc]{$Q\;=$}}
\multiput(13,3)(0,1){2}{\line(1,0){4}}
\put(13,2){\line(1,0){3}}
\multiput(13,0)(1,0){2}{\line(0,1){4}}
\multiput(13,0)(0,1){2}{\line(1,0){1}}
\multiput(15,2)(1,0){2}{\line(0,1){2}}
\put(17,3){\line(0,1){1}}
\put(3.5,3.5){\makebox(0,0)[cc]{\footnotesize 1}}
\put(4.5,3.5){\makebox(0,0)[cc]{\footnotesize 3}}
\put(5.5,3.5){\makebox(0,0)[cc]{\footnotesize 6}}
\put(6.5,3.5){\makebox(0,0)[cc]{\footnotesize 7}}
\put(3.5,2.5){\makebox(0,0)[cc]{\footnotesize 2}}
\put(4.5,2.5){\makebox(0,0)[cc]{\footnotesize 4}}
\put(5.5,2.5){\makebox(0,0)[cc]{\footnotesize 8}}
\put(3.5,1.5){\makebox(0,0)[cc]{\footnotesize 5}}
\put(3.5,0.5){\makebox(0,0)[cc]{\footnotesize 9}}
\put(13.5,3.5){\makebox(0,0)[cc]{\footnotesize 1}}
\put(14.5,3.5){\makebox(0,0)[cc]{\footnotesize 2}}
\put(15.5,3.5){\makebox(0,0)[cc]{\footnotesize 5}}
\put(16.5,3.5){\makebox(0,0)[cc]{\footnotesize 7}}
\put(13.5,2.5){\makebox(0,0)[cc]{\footnotesize 3}}
\put(14.5,2.5){\makebox(0,0)[cc]{\footnotesize 4}}
\put(15.5,2.5){\makebox(0,0)[cc]{\footnotesize 9}}
\put(13.5,1.5){\makebox(0,0)[cc]{\footnotesize 6}}
\put(13.5,0.5){\makebox(0,0)[cc]{\footnotesize 8}}
\end{picture}
\end{center}
of tableaux. We have $\inv(P)=8$ and $\inv(Q)=7$. The number of squares 
laying in the second or fourth row of $P$ equals 4. Thus $\pi$ is an odd permutation.
\erm
\end{exa}

\begin{rem}
\brm
In \cite{reifegerste}, the problem of how to obtain the sign from the 
associated tableaux was solved for $321$-avoiding permutations. 
Note that the tableaux have at most two rows in this case. Let $\pi\in{\cal S}_n$ 
be a $321$-avoiding permutation and $(P,Q)$ its associated pair of tableaux. 
\cite[Proposition 2.1]{reifegerste} states that $\sg(\pi)=(-1)^{s+e}$ where $s$ denotes the sum of all 
entries of the second row of $P$ and $Q$ and $e$ is the length of this row.\\ 
In general, $s$ is not equal to $\inv(P)+\inv(Q)$. For example, 
$\pi=2\:5\:1\:6\:8\:3\:9\:4\:7\in{\cal S}_9$ is in bijection with the pair
\begin{center}
% == Picture ==========================================
\unitlength=0.4cm
\begin{picture}(19,2)
\linethickness{0.45pt}
\put(0,1){\makebox(0,0)[lc]{$P\;=$}}
\multiput(3,1)(0,1){2}{\line(1,0){5}}
\put(3,0){\line(1,0){4}}
\multiput(3,0)(1,0){5}{\line(0,1){2}}
\put(8,1){\line(0,1){1}}
\put(11,1){\makebox(0,0)[lc]{$Q\;=$}}
\multiput(14,1)(0,1){2}{\line(1,0){5}}
\put(14,0){\line(1,0){4}}
\multiput(14,0)(1,0){5}{\line(0,1){2}}
\put(19,1){\line(0,1){1}}
\put(3.5,1.5){\makebox(0,0)[cc]{\footnotesize 1}}
\put(4.5,1.5){\makebox(0,0)[cc]{\footnotesize 3}}
\put(5.5,1.5){\makebox(0,0)[cc]{\footnotesize 4}}
\put(6.5,1.5){\makebox(0,0)[cc]{\footnotesize 7}}
\put(7.5,1.5){\makebox(0,0)[cc]{\footnotesize 9}}
\put(3.5,0.5){\makebox(0,0)[cc]{\footnotesize 2}}
\put(4.5,0.5){\makebox(0,0)[cc]{\footnotesize 5}}
\put(5.5,0.5){\makebox(0,0)[cc]{\footnotesize 6}}
\put(6.5,0.5){\makebox(0,0)[cc]{\footnotesize 8}}
\put(14.5,1.5){\makebox(0,0)[cc]{\footnotesize 1}}
\put(15.5,1.5){\makebox(0,0)[cc]{\footnotesize 2}}
\put(16.5,1.5){\makebox(0,0)[cc]{\footnotesize 4}}
\put(17.5,1.5){\makebox(0,0)[cc]{\footnotesize 5}}
\put(18.5,1.5){\makebox(0,0)[cc]{\footnotesize 7}}
\put(14.5,0.5){\makebox(0,0)[cc]{\footnotesize 3}}
\put(15.5,0.5){\makebox(0,0)[cc]{\footnotesize 6}}
\put(16.5,0.5){\makebox(0,0)[cc]{\footnotesize 8}}
\put(17.5,0.5){\makebox(0,0)[cc]{\footnotesize 9}}
\end{picture}
\end{center}
of tableaux for which $s=47$. On the other hand, $\inv(P)=9$ and $\inv(Q)=4$.\\
Alternatively to the proof given in \cite{reifegerste}, we can derive the 
result from Theorem \ref{main}. Let $i_k$ be the entry occupying in the $k$th square of the first row of $P$. 
Then $i_k-k$ of the elements contained in the second row are smaller than 
$i_k$. Similarly, the entry $j_k$ of the $k$th square in the first row of 
$Q$ causes $j_k-k$ inversions. Since 
$s=n(n+1)-(i_1+\ldots+i_{n-e}+j_1+\ldots+j_{n-e})$, we have
\bdpm
\inv(P)+\inv(Q)=\sum_{k=1}^{n-e} (i_k-k)+\sum_{k=1}^{n-e} (j_k-k)\equiv s\bmod 2. 
\edpm  
\erm
\end{rem}
\vspace*{0.5cm}

% =======================================================================
% Section 5: Consequences
% =======================================================================

\setcounter{section}{5}\setcounter{theo}{0}

\centerline{\large{\bf 5}\hspace*{0.25cm}
{\sc Consequences}}
\vspace*{0.5cm}

Recently, Stanley \cite{stanley1} considered sign-balanced posets and gave 
particularly some new results for posets arising from partitions. (For an 
introduction to the theory of posets see \cite{stanley}.)\\
For a partition $\lambda$, let
\bdpm
I_\lambda(q)=\sum_T q^{\inv(T)} 
\edpm
be the generating function for the tableaux of shape $\lambda$ by the number of inversions. 
The integer $I_\lambda=I_\lambda(-1)$ is called the {\it imbalance} of $\lambda$.\\[2ex]
In this section, we discuss some consequences of the relation between the sign 
of a permutation and the signs of its associated tableaux. On the one hand, we 
give a simple proof for one of Stanley's results on imbalances of 
partitions. On the other hand, we can use the knowledge about 
imbalances of special shapes for refining the enumeration of certain restricted 
permutations concerning their sign and the length of their longest increasing 
subsequence.\\[2ex]
There are as many even permutations in ${\cal S}_n$ as odd ones. Applying Theorem 
\ref{main}, we can interpret this well-known fact in terms of tableaux. This 
extends one of Stanley's results and proves a special case of a conjecture he 
has given. By \cite[Theorem 3.2b]{stanley1}, we have 
\bdpm
\sum_{\lambda\vdash 2m} (-1)^{v(\lambda)} I_\lambda^2=0
\edpm 
where $v(\lambda)$ denotes the maximum number of disjoint vertical dominos that 
fit in the shape $\lambda$. The proof uses a bijection between colored 
biwords and pairs of standard domino tableaux of the same shape which was 
established by Shimozono and White \cite[Theorem 30]{shimozono-white}.\\
Note that $v(\lambda)$ just counts the number of squares contained in an 
even-indexed row of $\lambda$:  
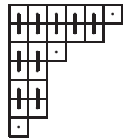
\begin{center}
% == Picture ==========================================
\unitlength=0.25cm
\begin{picture}(6,7)
\multiput(0,0)(1,0){2}{\line(0,1){7}}
\put(2,1){\line(0,1){6}}
\put(3,4){\line(0,1){3}}
\multiput(4,5)(1,0){2}{\line(0,1){2}}
\put(6,6){\line(0,1){1}}
\multiput(0,6)(0,1){2}{\line(1,0){6}}
\put(0,5){\line(1,0){5}}
\put(0,4){\line(1,0){3}}
\multiput(0,1)(0,1){3}{\line(1,0){2}}
\put(0,0){\line(1,0){1}}
\multiput(0.5,6.5)(1,0){6}{\circle*{0.2}}
\multiput(0.5,5.5)(1,0){5}{\circle*{0.2}}
\multiput(0.5,4.5)(1,0){3}{\circle*{0.2}}
\multiput(0.5,3.5)(1,0){2}{\circle*{0.2}}
\multiput(0.5,2.5)(1,0){2}{\circle*{0.2}}
\multiput(0.5,1.5)(1,0){2}{\circle*{0.2}}
\put(0.5,0.5){\circle*{0.2}}
\linethickness{1pt}
\multiput(0.5,5.5)(1,0){5}{\line(0,1){1}}
\multiput(0.5,3.5)(1,0){2}{\line(0,1){1}}
\multiput(0.5,1.5)(1,0){2}{\line(0,1){1}}
\end{picture}
\end{center}  
\vspace*{-1.5ex}

To emphasize the equivalent definition, we rename the statistic $v$ and write 
$e(\lambda)$ to denote the sum $\lambda_2+\lambda_4+\ldots$ of all even-indexed parts 
of $\lambda$.\\[2ex]
The case $t=1$ of Conjecture \cite[3.3b]{stanley1} claims that the above equation is correct for partitions of an arbitrary integer $n$. Theorem 
\ref{main} yields the proof.

\begin{theo}
For all $n\ge 1$, we have
\bdpm
\sum_{\lambda\vdash n} (-1)^{e(\lambda)} I_\lambda^2=0.
\edpm
\end{theo}

\begin{proof} From the sign-balance on ${\cal S}_n$ we obtain
\bdpm
0=\sum_{\pi\in{\cal S}_n} \sg(\pi)=\sum_{\lambda\vdash n}\sum_{(P,Q)} (-1)^{e(\lambda)}\sg(P)\sg(Q)=
\sum_{\lambda\vdash n} (-1)^{e(\lambda)} I_\lambda^2 
\edpm
where $(P,Q)$ ranges over all pairs of tableaux of shape $\lambda$.
\end{proof}

For certain shapes $\lambda$ (hooks and rectangles), the imbalance 
$I_\lambda$ has been determined explicitly. The characterization of permutations whose associated tableaux have exactly 
such a shape allows a refined 
enumeration of these permutations regarding their sign 
and even the length of their longest increasing subsequence. These 
considerations lead to pattern-avoiding permutations. Given a permutation $\pi\in{\cal S}_n$ and a permutation $\tau\in{\cal S}_k$, 
we say that $\pi$ {\it avoids the pattern} $\tau$ if there is no sequence $1\le 
i_1<i_2<\ldots<i_k\le n$ such that the elements 
$\pi_{i_1},\pi_{i_2},\ldots,\pi_{i_k}$ are in the same relative order as 
$\tau_1,\tau_2,\ldots,\tau_k$.\\
In \cite[Theorem 1.1]{reifegerste}, the joint distribution of $\sg$ and $\lis$, 
the length of the longest increasing subsequence, was given for 
$321$-avoiding permutations. As an application of result \cite[Proposition 3.4]{stanley1} which deals with the weighted imbalances of hooks, we obtain the number of even 
and odd permutations which avoid both $213$ and $231$ respecting the 
statistic $\lis$, respectively.  

\begin{prop} \label{hooks}
Let $\pi\in{\cal S}_n$ be a permutation and $\lambda$ the shape of its 
insertion tableau $P$. Then $\pi$ avoids $213$ and $231$ if and only if 
$w(P)=12\cdots n$ and $\lambda=(k,1^{n-k})$ for $k\in\{1,\ldots,n\}$. 
\end{prop}

\begin{proof}
First we show that permutations having insertion tableaux of this very special form can be characterized as 
follows. Let $i_1<i_2<\ldots<i_d$ be the descents of $\pi\in{\cal S}_n$ and 
$j_1<j_2<\ldots<j_{n-d}$ the remaining positions. (An integer 
$i\in\{1,\ldots,n-1\}$ is called a {\it descent} of $\pi$ if $\pi_i>\pi_{i+1}$.)  
Then the insertion tableau $P$ of $\pi$ has no inversions and is of shape $(k,1^{n-k})$
if and only if $d=n-k$ and
\bdpm
\pi_{i_1}\pi_{i_2}\cdots\pi_{i_d}=n(n-1)\cdots(k+1)\quad\mbox{and}\quad
\pi_{j_1}\pi_{j_2}\cdots\pi_{j_{n-d}}=12\cdots k.
\edpm
Suppose that 
\begin{center}
% == Picture ==========================================
\unitlength=0.35cm
\begin{picture}(9,5)
\put(0,2.5){\makebox(0,0)[lc]{$P\;=$}}
\multiput(3,0)(1,0){2}{\line(0,1){5}}
\multiput(8,4)(1,0){2}{\line(0,1){1}}
\put(5,4){\line(0,1){1}}
\multiput(3,4)(0,1){2}{\line(1,0){6}}
\put(3,3){\line(1,0){1}}
\multiput(3,0)(0,1){2}{\line(1,0){1}}
\put(6.5,4.5){\makebox(0,0)[cc]{$\cdots$}}
\put(3.5,2.25){\makebox(0,0)[cc]{$\vdots$}}
\put(3.5,4.5){\makebox(0,0)[cc]{\footnotesize 1}}
\put(4.5,4.5){\makebox(0,0)[cc]{\footnotesize 2}}
\put(8.5,4.5){\makebox(0,0)[cc]{\footnotesize $k$}}
\put(3.5,0.5){\makebox(0,0)[cc]{\footnotesize $n$}}
\put(3.5,3.5){\makebox(0,0)[cc]{\footnotesize $l$}}
\end{picture}
\end{center}
where $l=k+1$. By reversing the steps in the Robinson-Schensted-Knuth algorithm, 
we obtain the elements $\pi_n,\pi_{n-1},\ldots,\pi_1$ from the entry contained 
in the rightmost square of the first row of the tableaux 
$P=P_n,P_{n-1},\ldots,P_1$. To see this note that the shape of the recording 
tableau $Q$ is a hook as well. Clearly, if the element $i$ appears at the end 
of the first row in the subtableau $Q_i$, then we find $\pi_i$ at the same 
position in $P_i$. Otherwise, if the element $i$ occupies the bottom square in the first 
column of $Q_i$, then by applying the reverse row-insertion to $P_i$ with that 
square, the rigtmost element is bumped out of the first row as well since all 
entries of the hook leg are larger than the entries of the hook arm. 
Therefore, a descent arises whenever an element is bumped out that was not placed 
in the first row of $P$. Since the reverse procedure builds up the permutation 
from right to left, the descent tops (that is, the elements 
$\pi_{i_1},\ldots,\pi_{i_{n-k}}$) are decreasing while the remaining letters form 
the sequence $1,2,\ldots,k$.\\
Conversely, it is obvious that a permutation with this descent structure 
produces an insertion tableau whose row word is the identity and which has 
the shape of a hook.\\
But these permutations are exactly the $\{213,231\}$-avoiding ones. If 
$\pi\in{\cal S}_n$ contains no pattern $213$ or $231$, then for any descent $i$ 
there exists no integer $j>i$ for which $\pi_j>\pi_i$ (otherwise we have the pattern 
$213$ in $\pi$), and for any non-descent $i$ there exists no $j>i$ with 
$\pi_j<\pi_i$ (otherwise the pattern $231$ occurs). Consequently, the descent 
tops of $\pi$ have to be decreasing, and moreover, equal to $n,n-1,\ldots,k+1$; 
the remaining elements are arranged in increasing order. By similar reasoning, 
we obtain the converse.
\end{proof}

\begin{cor} \label{refine}
Let $A_n$ denote the set of permutations in ${\cal S}_n$ which avoid both $213$ and 
$231$. For all $n\ge1$, we have
\bdpm
\sum_{\pi\in A_n}\sg(\pi)q^{\lis(\pi)}=q(q+1)^{\lf (n-1)/2\rf}(q-1)^{\lf n/2\rf}.
\edpm  
In particular, there are as many even as odd permutations in $A_n$. 
\end{cor}

\begin{proof}
Let $\pi\in A_n$ with $\lis(\pi)=k$. By Proposition \ref{hooks}, the tableaux in bijection with 
$\pi$ are of shape $\lambda=(k,1^{n-k})$. (Recall that the length of the longest increasing subsequence 
equals the length of the first row of the associated tableaux.) Clearly, 
$e(\lambda)=\lceil\frac{n-k}{2}\rceil$ and $e(\lambda')=\lf\frac{k}{2}\rf$ for the conjugate 
partition $\lambda'$. By the proof of \cite[Proposition 4.3]{stanley1}, we have
$I_\lambda=0$ if $n$ is odd and $k$ even and 
$I_\lambda={h(n)\choose h(k)}$ otherwise where 
$h(x)=\lf(x-1)/2\rf$. By means of Theorem \ref{main}, we can express the coefficients of the generating function 
\bdpm
F(q)=\sum_{\pi\in A_n} \sg(\pi)q^{\lis(\pi)}
\edpm
in terms of imbalances:
\bdpm
F(q)=\sum_{k=1}^n\sum_{Q} (-1)^{e(\lambda)}\sg(Q)q^k
=\sum_{k\,\mbox{\tiny odd}} q^{2e(\lambda')+1}(-1)^{e(\lambda)}I_\lambda+
\sum_{k\,\mbox{\tiny even}} q^{2e(\lambda')}(-1)^{e(\lambda)}I_\lambda
\edpm 
where $Q$ ranges over all tableaux of shape $\lambda=(k,1^{n-k})$. If $n$ is odd, then $I_\lambda$ vanishs whenever $k$ is even. Thus 
\bdpm
F(q)=q\sum_{k=1}^n q^{2e(\lambda')}(-1)^{e(\lambda)}I_\lambda=q(q^2-1)^{\frac{n-1}{2}}. 
\edpm 
In case of even $n$, we have
\bdpm
\sum_{k\,\mbox{\tiny even}} q^{2e(\lambda')}(-1)^{e(\lambda)}I_\lambda=
\sum_{k\,\mbox{\tiny even}} q^k(-1)^{\frac{n-k}{2}}
{\frac{n-2}{2}\choose\frac{k-2}{2}}=
-q\sum_{k\,\mbox{\tiny odd}} q^{2e(\lambda')+1}(-1)^{e(\lambda)}I_\lambda
\edpm
and hence
\bdpm
F(q)=q(1-q)\sum_{k\,\mbox{\tiny odd}} q^{2e(\lambda')}(-1)^{e(\lambda)}I_\lambda=
q(q-1)(q^2-1)^{\frac{n}{2}-1}.
\edpm 
\end{proof}

\begin{rem}
\brm
Proposition \ref{hooks} gives a bijection between $\{213,231\}$-avoiding 
permutations and standard Young tableaux of hook shape (consider only the recording 
tableaux). Therefore, it also yields a proof -- although not the most obvious one -- 
for the known fact $|A_n|=2^{n-1}$.\\
Because of the symmetry of the Robinson-Schensted-Knuth correspondence, we may 
replace the pattern $231$ in Corollary \ref{refine} with $312$.
\erm
\end{rem}
\vspace*{0.5cm}

{\bf Acknowledgement} I would like to thank Richard Stanley for helpful 
communications. 
\vspace*{0.75cm}

% =======================================================================
% References
% =======================================================================

\centerline{\large\sc References}
\vspace*{0.5cm}

\begin{enumbib}

\bibitem{beissinger}
J. S. Beissinger, 
Similar constructions for Young tableaux and involutions, and their application 
to shiftable tableaux, 
{\it Discrete Math.} {\bf 67} (1987), 149-163.

\bibitem{fulton}
W. Fulton, 
{\it Young tableaux}, 
Cambridge University Press, 1997.

\bibitem{knuth}
D. E. Knuth, 
Permutations, matrices, and generalized Young tableaux, 
{\it Pacific J. Math.} {\bf 34} (1970), 709-727.

\bibitem{reifegerste}
A. Reifegerste,
Refined sign-balance on 321-avoiding permutations, preprint, 2003, 
math.CO/0503327.

\bibitem{robinson}
G. de B. Robinson, 
On the representations of the symmetric group, 
{\it Amer. J. Math.} {\bf 60} (1938), 745-760.

\bibitem{schensted}
C. Schensted, 
Longest increasing and decreasing subsequences, 
{\it Canad. J. Math.} {\bf 13} (1961), 179-191.

\bibitem{schuetzenberger}
M. P. Sch\"utzenberger, 
Quelques remarques sur une construction de Schensted, 
{\it Math. Scand.} {\bf 12} (1963), 117-128.

\bibitem{schuetzenberger1}
M. P. Sch\"utzenberger, 
La correspondence de Robinson,
in {\it Combinatoire et repr\'esentation du groupe sym\'etrique}, Lecture Notes 
in Math. {\bf 579}, Springer, 1977, 59-113.

\bibitem{shimozono-white}
M. Shimozono and D. E. White, 
A Color-to-Spin Domino Schensted Algorithm, 
{\it Electron. J. Combinatorics} {\bf 8} (2001), R21. 

\bibitem{stanley}
R. P. Stanley, 
{\it Enumerative Combinatorics}, vol. I, II,
Cambridge University Press, 1997, 1999.

\bibitem{stanley1}
R. P. Stanley, 
Some Remarks on Sign-Balanced and Maj-Balanced Posets, 2002, math.CO/0211113.

\bibitem{viennot}
G. Viennot, 
Une forme g\'eom\'etrique de la correspondance de Robinson-Schensted, 
in {\it Combinatoire et repr\'esentation du groupe sym\'etrique}, Lecture Notes 
in Math. {\bf 579}, Springer, 1977, 29-58.

\end{enumbib}

\end{document}